\documentclass[a4paper,12pt]{article}
\usepackage{fullpage}
\usepackage[english,francais]{babel}
\usepackage{epsf}
\usepackage[dvips]{epsfig}
\usepackage{amssymb}
\usepackage{amsmath}
\usepackage[all]{xy}

\usepackage{ucs}
\usepackage[utf8x]{inputenc}\usepackage{wasysym}
\usepackage{aeguill}

\newtheorem{theorem}{Theorem}
\newtheorem{prop-f}[theoreme]{Proposition}

\newcommand{\finpreuve}{\hspace{\stretch{1}}{$\square$}}

\newcommand{\card}{{\rm card}}

\newcommand{\R}{\mathbb{R}}

\newcommand{\Z}{\mathbb{Z}}
\newcommand{\N}{\mathbb{N}}

\renewcommand{\epsilon}{\varepsilon}
\def\btab{\begin{eqnarray*}}
\def\etab{\end{eqnarray*}}
\def\beq{\begin{equation}}
\def\eeq{\end{equation}}

\newcounter{numeroexo}

\begin{document}

\selectlanguage{english}

\title{Monotonicity in first-passage percolation}

\author{
Jean-Baptiste Gou\'er\'e
\footnote{
\noindent\textit{Postal address}:
Universit\'e d'Orl\'eans
MAPMO
B.P. 6759
45067 Orl\'eans Cedex 2
France
\textit{E-mail}: jbgouere@univ-orleans.fr
}
}


\date{}
\maketitle

\begin{abstract}
We consider standard first-passage percolation on $\Z^d$.
Let $e_1$ be the first coordinate vector.
Let $a(n)$ be the expected passage time from the origin to $ne_1$.
In this short paper, we note that $a(n)$ is increasing under some strong condition on the support of the distribution of the passage times on the edges. 
\end{abstract}

\section{Introduction and results}

\paragraph*{First passage percolation.}

We consider the graph $\Z^d$, $d\ge 2$, obtained by taking $\Z^d$ as vertex set 
and by puting an edge between two vertices if the Euclidean distance between them is $1$.
We consider a family of non-negative i.i.d.r.v.\ $\tau={(\tau(e))}_{e \in {\cal E}}$ indexed by the set of edges ${\cal E}$ of the graph.
We interpret $\tau(e)$ as the time needed to travel along the edge $e$ (the graph is unoriented).

If $a$ and $b$ are two vertices of $\Z^d$, we call path from $a$ to $b$ any finite sequence of vertices $r=(a=x_0,...,x_k=b)$ such that,
for all $i\in \{0,...,k-1\}$, the vertices $x_i$ et $x_{i+1}$ are linked by an edge.
We denote by ${\cal C}(a,b)$ the set of such paths.
The time needed to travel along a path $r=(x_0,...,x_k)$ is defined by:
$$
\tau(r)=\sum_{i=0}^{k-1} \tau(x_i,x_{i+1}).
$$
Then, the time needed to go from $a$ to $b$ is defined by:
$$
T(a,b)=\inf\{\tau(r) : r \in {\cal C}(a,b)\}.
$$

Let $e_1, \dots, e_d$ denote the canonical basis vectors of $\R^d$.
We are interested in the sequence $(a(n))$ defined by :
$$
a(n)=E(T(0,ne_1)).
$$
We write $T'(0,ne_1)$ and $a'(n)$ for the passage times and expected passage times obtained when the paths are restricted to $\{(x_1,\dots,x_d) : 0 \le x_1 \le n\}$.

\paragraph*{Main result and related results.}

We denote by $S_-$ the infimum of the support of the distribution of the $\tau(e)$.
We denote by $S_+$ the supremum of the support.

\begin{theorem} \label{t} Assume $0<S_-$ and $S_+ \le 2S_-$. Then the sequence $(a_n)$ is non-decreasing. More precisely, we have:
$$
a(n) \ge a(n-1)+S_-\left[1-\frac{(S_+-S_-)^2}{{S_-}^2}\right].
$$
\end{theorem}

As soon as the distribution of the $\tau(e)$ is not a Dirac distribution there exists,
with probability one, infinitely many random $N$ such that
\footnote{Let us sketch a proof.
Fix $a$ and $b$ such that $S_-<a<b<S_+$. 
For each $n$, consider a box $\{n-C,\dots,n\} \times \{-D,\dots,D\}^{d-1}$.
Let $A_n$ be the following event: $\tau(e) \le a$ for edges inside the boundary of the box and $\tau(e) \ge b$ for edges inside the box.
For suitably chosen large $C$ and $D$ and for $n>C$, we have $T(0,(n-1)e_1) > T(0,ne_1)$ as soon as $A_n$ occurs.
As the $A_n$ are local event of fixed positive probability, the result follows.}
$T(0,(N-1)e_1)>T(0,Ne_1)$.
However, monotonicity of expected passage times seems quite natural and was already conjectured by Hammersley and Welsh in \cite{Hammersley-Welsh}.	 
In \cite{Alm-Wierman}, Alm and Wierman proved the monotonicity for $\Z \times \N$ and other $2$ dimensional models.
In \cite{Ahlberg-fppdim1}, Ahlberg made a detailed study of first passage percolation on essentially one-dimensional graphs,
an example of which is $\Z \times \{0,\dots,K\}^{d-1}$.
In particular, he proved the existence of a constant $n_0$, depending on the graph, such that $n \ge n_0$ implies $a(n) \ge a(n-1)$.
In \cite{Howard-monotone}, Howard proved the monotocity for an Euclidean first-passage percolation model.
We are not aware of any other positive results.

On the other hand, van den Berg proved in \cite{vdBerg-nonmonotone} that, when $d=2$, one has $a'(2) < a'(1)$ when
$\tau(e)=1$ with small probability and $\tau(e)=0$ otherwise.
Note that we still have $a'(2)<a'(1)$ if, instead of setting $\tau(e)=0$ we set $\tau(e)=\epsilon$ for a small enough
\footnote{Indeed, $a'(1)$ can only increase while $a'(2)$ increases by at most $\epsilon E(N)$ where $N$ is the length of a geodesic for the initial passage times.
Using $T'(0,2e_1) \le 2$ one can check that any geodesic must remain in a random box of subgeometrical height.
Therefore $E(N)$ is finite and the result follows.} $\epsilon$.
A related result was given by Joshi in \cite{Joshi}.

We refer to the review by Howard \cite{Howard-revue-fpp} for a more detailed account.

\paragraph*{Further remarks.}

\begin{itemize}
\item The same result holds for the $a'(n)$.
\item The proof gives that $T(0,ne_1)$ stochastically dominates the mean of $n$ dependent copies of $T(0,(n-1)e_1)$
(see \eqref{e:domination} and \eqref{e:memeloi}).
\item With the same strategy one can prove for example the following result:
\begin{equation}\label{e:CS2}
a(n) \ge a(n-1) \hbox{ as soon as } S_->0 \hbox{ and } \left(\frac{a(n)n^{-1}}{S_-}-1\right)\left(\frac{E(\tau(e))}{S_-}-1\right) \le \frac12
\end{equation}
where $e$ is a fixed edge. 
We show how to adapt the proof of Theorem \ref{t} to prove this result below the proof of Theorem \ref{t}.
In particular, using the inequality $a(n) \le nE(\tau(e))$, we get that $a$ is non-decreasing as soon as :
$$
S_->0 \hbox{ and } E(\tau(e)) \le (1+2^{-1/2})S_-.
$$
This gives a sufficient condition with no assumption on $S_+$ which can be infinite.
However, this sufficient condition is still strong and we do not see how to give any significantly weaker condition.
\item Fix the distribution of $\tau(e)$. Assume $S_->0$ and $E\tau(e)<\infty$. 
Then the conditions in \eqref{e:CS2} are true for large enough $n$ and $d$.
This is due to the fact that $a(n)n^{-1}$ can be made arbitrarily close to $S_-$.
\end{itemize}

\section{Proofs}

\paragraph{Proof of Theorem \ref{t}.}
 For all $i$ we consider the following sets of edges:
\begin{itemize}
\item $H^i$: the set of edges $(x,x+e_1)$ where $x=(x_1,\dots,x_d)$ is such that $x_1=i$.
\item $V^i$: the set of edges $(x,x+e_k)$ where $x_1=i$ and $k$ belongs to $\{2,\dots,d\}$.
\end{itemize}
We define new passage times $\tau^i(e)$ as follows:
\begin{itemize}
\item If $e$ belongs to $H^i$ then $\tau^i(e)=0$.
\item If $e$ belongs to $V^i$ then $\tau^i(e)=+\infty$.
\item Otherwise, $\tau^i(e)=\tau(e)$.
\end{itemize}
We denote by $T^i(a,b)$ the time needed to travel from $a$ to $b$ with the passage times $\tau^i(e)$.
Note, for all $n\ge 1$ and all $i \in \{0,n-1\}$, the following:
\begin{equation}\label{e:memeloi}
T^i(0,ne_1) \hbox{ and } T(0,(n-1)e_1) \hbox{ have the same distribution.}
\end{equation}

We now compare $T^i(0,ne_1)$ and $T(0,ne_1)$.
Let $\pi$ be a path from $0$ to $ne_1$ such that $\tau(\pi)=T(0,ne_1)$.
We modify this path as follows.
Each time the path goes, in this order, through an edge $(x,y)\in V^i$, 
we replace this part of the path by $(x,x+e_1,y+e_1,y)$.
We denote by $\pi^i$ the modified path.
We have
$$
\tau^i(\pi^i) \le \tau(\pi) - S_- \card(\pi \cap H_i) + (S_+-S_-) \card(\pi \cap V_i) 
$$
where, for example, $\card(\pi \cap H_i)$ denotes the number of edges of $H_i$ used by $\pi$.
The term involving $H_i$ is due to the time saved by the modification of the passage times.
The term involving $V_i$ is partly due to the time left by the modification of the path.
We thus get
\begin{equation}\label{e:ti}
T^i(0,ne_1) \le T(0,ne_1) - S_- \card(\pi \cap H_i) + (S_+-S_-) \card(\pi \cap V_i)
\end{equation}
and then
\begin{equation}\label{e:somme}
\sum_{i=0}^{n-1} T^i(0,ne_1) \le nT(0,ne_1) - S_- \sum_{i=0}^{n-1}\card(\pi \cap H_i) + (S_+-S_-) \sum_{i=0}^{n-1}\card(\pi \cap V_i).
\end{equation}
Note 
\begin{equation}\label{e:inegaliteH}
\sum_{i=0}^{n-1}\card(\pi \cap H_i) \ge n,
\end{equation}
as $\pi$ is a path from $0$ to $ne_1$.
But
\begin{eqnarray}
T(0,ne_1) 
 & = & \tau(\pi) \nonumber \\
 & \ge & S_- \sum_{i=0}^{n-1}\card(\pi \cap V_i)+ S_-\sum_{i=0}^{n-1}\card(\pi \cap H_i)  \nonumber \\
 & \ge & S_- \sum_{i=0}^{n-1}\card(\pi \cap V_i)+ S_-n \label{e:inegaliteV}
\end{eqnarray}
and, moreover,
\begin{eqnarray*}
T(0,ne_1) 
 & \le & \tau(0,e_1,2e_2,\dots,ne_1) \\
 & \le & n S_+.
\end{eqnarray*}
Therefore:
\begin{eqnarray}
\sum_{i=0}^{n-1}\card(\pi \cap V_i) 
 & \le & \frac{T(0,ne_1) - nS_-}{S_-} \label{e:majorationV} \\
 & \le & \frac{nS_+-nS_-}{S_-}. \nonumber
\end{eqnarray}
From \eqref{e:somme} and \eqref{e:majorationV} we get:
\begin{equation}\label{e:domination}
\sum_{i=0}^{n-1} T^i(0,ne_1) \le nT(0,ne_1) - nS_-  + \frac{n(S_+-S_-)^2}{S_-}.
\end{equation}
Taking expectations and using \eqref{e:memeloi} we get:
$$
n a(n-1) \le n a(n)	 -n S_-\left[1-\frac{(S_+-S_-)^2}{{S_-}^2}\right].
$$
The proof follows. \finpreuve

\paragraph{Proof of \eqref{e:CS2}.} 
The proof is essentially the same.
The main difference lies in the definition of the new passage times $\tau^i(e)$.
We let $\widetilde{\tau}$ be an independent copy of $\tau$.
We then set:
\begin{itemize}
\item If $e$ belongs to $H^i$ then $\tau^i(e)=0$.
\item If $e$ belongs to $V^i$ then $\tau^i(e)=+\infty$.
\item If $e$ belongs to $V^{i+1}$ then $\tau^i(e)=\widetilde\tau^i(e)$.
\item Otherwise, $\tau^i(e)=\tau(e)$.
\end{itemize}
Instead of \eqref{e:ti} we can write, after taking conditional expectation w.r.t.\ $\tau$ :
\begin{equation}\label{e:ti2}
T^i(0,ne_1) \le T(0,ne_1)  - S_- \card(\pi \cap H_i) + (E(\tau(e))-S_-) \card(\pi \cap V_i) + (E(\tau(e))-S_-) \card(\pi \cap V_{i+1}).
\end{equation}
Instead of \eqref{e:somme} we can get :
\begin{equation}\label{e:somme2}
\sum_{i=0}^{n-1} T^i(0,ne_1) \le nT(0,ne_1) - S_- \sum_{i=0}^{n-1}\card(\pi \cap H_i) + 2(E(\tau(e))-S_-) \sum_{i=0}^n\card(\pi \cap V_i).
\end{equation}	
Using \eqref{e:inegaliteH}, an equality similar to \eqref{e:inegaliteV} and taking expectation, we get:
$$
na(n-1) \le na(n) - S_- n + 2(E(\tau(e))-S_-) \frac{a(n)-nS_-}{S_-}
$$
and thus:
$$
a(n-1) \le a(n) - S_- \left(1-2\left(\frac{E(\tau(e))}{S_-}-1\right)\left(\frac{a(n)n^{-1}}{S_-}-1\right)\right).
$$
The proof follows. \finpreuve

\def\cprime{$'$} \def\cprime{$'$}


\begin{thebibliography}{1}

\bibitem{Ahlberg-fppdim1}
Daniel Ahlberg.
\newblock Asymptotics of first-passage percolation on 1-dimensional graph.
\newblock {\em arXiv:1107.2276}, 2011.

\bibitem{Alm-Wierman}
Sven~Erick Alm and John~C. Wierman.
\newblock Inequalities for means of restricted first-passage times in
  percolation theory.
\newblock {\em Combin. Probab. Comput.}, 8(4):307--315, 1999.
\newblock Random graphs and combinatorial structures (Oberwolfach, 1997).

\bibitem{Hammersley-Welsh}
J.~M. Hammersley and D.~J.~A. Welsh.
\newblock First-passage percolation, subadditive processes, stochastic
  networks, and generalized renewal theory.
\newblock In {\em Proc. {I}nternat. {R}es. {S}emin., {S}tatist. {L}ab., {U}niv.
  {C}alifornia, {B}erkeley, {C}alif}, pages 61--110. Springer-Verlag, New York,
  1965.

\bibitem{Howard-monotone}
C.~Douglas Howard.
\newblock Differentiability and monotonicity of expected passage time in
  {E}uclidean first-passage percolation.
\newblock {\em J. Appl. Probab.}, 38(4):815--827, 2001.

\bibitem{Howard-revue-fpp}
C.~Douglas Howard.
\newblock Models of first-passage percolation.
\newblock In {\em Probability on discrete structures}, volume 110 of {\em
  Encyclopaedia Math. Sci.}, pages 125--173. Springer, Berlin, 2004.

\bibitem{Joshi}
V.~M. Joshi.
\newblock First-passage percolation on the plane square lattice.
\newblock {\em Sankhy\=a Ser. A}, 39(2):206--209, 1977.

\bibitem{vdBerg-nonmonotone}
J.~van~den Berg.
\newblock A counterexample to a conjecture of {J}. {M}. {H}ammersley and {D}.
  {J}. {A}. {W}elsh concerning first-passage percolation.
\newblock {\em Adv. in Appl. Probab.}, 15(2):465--467, 1983.

\end{thebibliography}
\end{document}